\documentstyle{amsppt}
\magnification=\magstep1
\vsize=22 true cm
\hsize=16 true cm
\voffset=-1 true cm
\pageno=1
\topmatter
\title
{Generalized Hopfian property,\\
 minimal Haken manifold, and\\
  J. Simon's conjecture for 3-manifold groups}
\endtitle
\author{Alan W. Reid and Shicheng Wang and Qing Zhou}
\endauthor
\thanks{The first author was supported by the Royal Society, the
NSF and The Alfred P. Sloan Foundation.
The second and third author are supported by Outstanding Youth Fellowship
of NSFC.}
\endthanks

\abstract {We address a conjecture that $\pi_1$-surjective maps between
closed aspherical 3-manifolds having the same rank on $\pi_1$
must be of non-zero degree.
The conjecture is proved for Seifert manifolds,
which is used
in constructing the first known example of minimum Haken manifold.
Another motivation is to study epimorphisms of 3-manifold
groups via maps of non-zero degree between 3-manifolds.}
\endabstract

\leftheadtext{Generalized Hopfian Property}
\rightheadtext{A.W. Reid, S.C.Wang and Q.Zhou}
\endtopmatter
\document


{\bf Section 1. Introduction and some examples.}

Let $M$ and $N$ be closed 3-manifolds and 
$f : M \to N$ a map of non-zero degree, then
the image of $f_*$ is a subgroup of finite index in $\pi_1(N)$. If $M$
and $N$ are aspherical, any homomorphism $\phi : \pi_1(M)\to \pi_1(N)$
determines a unique map
$f : M\to N$ up to homotopy such that $f_*=\phi$.  It seems natural to
ask when there exists $f:M\to N$ of non-zero degree given
a homomorphism $\phi$ surjecting $\pi_1 (M)$ on a subgroup of finite
index in $\pi_1(N)$?  There are elementary
constructions of examples (see below) that show in general that the
answer is no.  Before discussing some examples we make the following
definition.

{\bf Definition 1.1.} A map $f:M\to N$ between 3-manifolds is
$\pi_1$-surjective (resp. $\pi_1$-finite-index) if $f_*: \pi_1(M)\to
\pi_1(N)$ is surjective (resp. the image of $f_*$ is a subgroup of
finite index).

Recall that if $M$ is an $n$-manifold, the { rank} of $\pi_1(M)$ (or simply
by abuse just $M$) is the minimal cardinality of a generating system for
$\pi_1(M)$.

Let us first have a look of the situation in dimension 2 which is quite simple.
Throughout the paper $\Sigma_k$ will denote a 
closed orientable surface of genus $k$.

{\bf Example 1.1.}
It is not difficult to see that there is a $\pi_1$-surjective map
$f: \Sigma_l\to \Sigma_k$ which is of degree zero when $l\ge 2k$.

On the otherhand, if $f: \Sigma_l\to \Sigma_k$ is a $\pi_1$-surjective
map with  $0<l< 2k$, then we claim that $f$ is of non-zero degree.
The proof of this result is direct.
Choose a 1-skeleton of $\Sigma_k$ to be a one point wedge of $2k$
circles ${\Cal V} = \vee C_i$. Fix a point $x_i$ on $C_i$.
If $f$ is of degree zero,  then the image of $f$ can be deformed into
$\Cal V$. We assume therefore that this is the case.
Since $f_*$ is $\pi_1$-surjective, $f: \Sigma_l\to \Cal V$ must be
surjective. We may also assume that $f$ is transverse to each $x_i$,
$i=1,2,..., 2k$. So $f^{-1}(\cup x_i)$ is a set of essential circles.
Partition $f^{-1}(\cup x_i)$ into sets $G_1, ...,G_h$ such that
two components are in the same set if and only if they are
parallel. For each $G_j$, find an annulus $A_j$ containing $G_j$.
Then squeeze each $A_j$ to an arc $a_j$ and the part 
$\Sigma_k \setminus\cup A_j$
to a point. The quotient $Q$ will be a bouquet of $h$ circles.
Since $\Cal V -\{x_i, i=1,2,..., 2k\}$ is contractible, the map
$f: \Sigma_l\to \Cal V$ factors through $q: Q\to \Cal V$ which is still
$\pi_1$-surjective.
It follows that $h\ge 2k$. In particular, there are at least $h\ge 2k$
disjoint essential non-separating non-parallel circles.
By a well known argument in surface topology, we must have that
the $l$, the genus of $\Sigma_l$, is at least $2k$.
$\square$

Let us come  back to dimension 3.
The first example illustrates the aspherical assumption.
\smallskip
\leftline{\bf \hbox{Example 1.2}:}
Let $f=e\circ p: S^2\times S^1\to S^2\times S^1$,
where $p$ is a map which pinches $S^2\times S^1$ to $ S^1$,
and $e$ identifies  $S^1$ to a fiber $*\times S^1\subset S^2\times S^1$.
Clearly $f$ is of zero degree but $\pi_1$-surjective.
$\square$

The second example shows that if we do not require that the
manifolds have the same rank, then the answer to the question is no.

\leftline{\bf \hbox{Example 1.3}:}
We construct a map $f : \Sigma_{g+1}\times S^1\to \Sigma_g\times S^1$ 
of zero degree which is $\pi_1$-surjective.
The map $f$ is the composition of the following four geometric operations.

\item{} Project $\Sigma_{g+1}\times S^1$ to $\Sigma_{g+1}$.

\item{} Squeeze a suitable separating circle on $\Sigma_{g+1}$
to a point in such a way that the quotient space is a one point union
of a torus and $\Sigma_g$.

\item{} Squeeze the torus to a circle
in such a way that the quotient space is a one point union
of the circle and $\Sigma_g$.

\item{} Send $\Sigma_g$ and the circle
to a section $\Sigma_g\times *$ and the circle fiber of 
$\Sigma_g\times S^1$ respectively.$\square$

The third example has the same purpose as the second one, but the manifolds
in this case are hyperbolic.

\leftline{\bf \hbox{Example 1.4}:}
Let $M$ be a closed hyperbolic 3-manifold whose fundamental group surjects
the free group of rank $2$. Such examples are easily constructed
by doing hyperbolic surgery on a null-homotopic hyperbolic knot
in $\Sigma_2\times S^1$ [Section 3, BW].
Let $\phi_1: \pi_1(M)\to {\Bbb Z}*{\Bbb Z}$ denote such a map.
Let $N$ be any hyperbolic 3-manifold such that $\pi_1(N)$ has
two generators, then there is an epimorphism  $\phi_2: {\Bbb Z}*{\Bbb Z}
\to \pi_1(N)$.
If we choose $N$ such that the volume of $N$ is larger than the volume
of $M$, then the map realizing the epimorphism $\phi=\phi_2\circ \phi_1$
must be zero degree by the work of Gromov and Thurston, \cite{T}.
We remark that the volumes of hyperbolic 3-manifolds of rank $2$
are unbounded. Briefly, it follows from work of Adams that the
volumes of hyperbolic 2-bridge knot complements are unbounded. Doing large
enough hyperbolic Dehn surgeries on these gives the required family, see [CR].

In fact it can be seen directly that the map realizing $\phi$ must be
of zero degree since such a map factors through a 1-dimensional
complex.  $\square$

As a consequence of these examples, we state the following more refined 
version of the question posed above.

\proclaim{\bf Question 1.5} Let $M$ and $N$ be closed aspherical
3-manifolds such that the rank of $\pi_1(M)$ equals the rank of
$\pi_1(N)$. Assume, that $\phi : \pi_1(M)\to \pi_1(N)$ is surjective
or whose image is a subgroup of finite index.  Does $\phi$ determine a
map $f : M\to N$ of non-zero degree?  \endproclaim
Note if $M$ and $N$
are homeomorphic and satisfy Thurston's geometrization conjecture,
then a $\pi_1$-surjective map $f: M\to M$ must be degree one.  For
since $\pi_1(M)$ is hopfian, $f_*$ is surjective implies $f_*$ is an
isomorphism.  Since $M$ is aspherical $f$ must be a homotopy
equivalence, and so in particular, $f$ is of degree one.  Thus the
question above is a kind of generalization of the Hopfian property:
the condition ``homeomorphic manifolds'' is replaced by ``manifolds of
the same rank'', the condition ``$\pi_1$-surjective'' is replaced by
``$\pi$-surjective'' or ``$\pi_1$-finite-index'', and the conclusion replace
``degree one'' by ``non-zero degree''.  It is easy to construct examples
to show that 
``non-zero degree'' cannot be sharpened to ``degree one'', see the examples
in Section 3.
\smallskip

One of the main results of this paper is to prove that for Seifert
fibered 3-manifolds Question 1.5 has a positive answer (see Theorem
2.1 and Remark 2.4).
In \S 4 we use this result to construct the
first known example of a Haken 3-manifold which is minimal with
respect to degree 1 mappings in Thurston's picture of 3-manifolds
(Theorem 4.1).
The manifold is a graph
manifold built from the union of two trefoil knot complements.
An orientable 3-manifold $M$ is minimal if there is a degree one map
$f: M\to N$ implies either $N=S^3$ or $M=N$. Usually it is difficult
tell if a 3-manifold is minimal.
We remark that all minimal Seifert manifolds are non-Haken \cite{LWZ},
and that the known minimal hyperbolic 3-manifolds are also non-Haken
\cite{RW}, see [BW], [RW] and [LWZ]
for a further discussion of such matters. 

We were also motivated by the following posed by J. Simon.

\proclaim {[K. Problem 1.12]}
Let $G_K=\pi_1(S^3-K)$ for a knot $K$ in $S^3$. Conjectures:
if there is an epimorphism $\phi : G_K\to G_L$, then

(A) rank $G_L>$ rank $G_K$.

(B) genus$(L)\ge$ genus$(K)$.

(C) Given $K$, there is a number $N_K$ such that any sequence of epimorphisms of
knot groups $G_K\to G_{L_1}\to .... \to G_{L_n}$ with $n\ge N_K$ contains
an isomorphism.

(D) Given $K$, there are are only finitely many knot groups $G$ for which
there is an epimorphism $G_K\to G$.
\endproclaim

These conjectures have seen little progress.  On the otherhand, more
recently, questions similar to (C) and (D) have been raised for degree one
maps and there are already several substantial results in this setting.

\proclaim{[K Problem 3.100]}
Let $M$ be a closed orientable 3-manifold.

(A) Are there only finitely many irreducible 3-manifolds $N$ such that
there exists a degree one map $M\to N$?

(B) Does there exist an integer $N_{M}$ such that if a sequence of degree one
map $M\to M_1\to ....\to M_k$ with $k>N_{M_0}$, the sequence contains an
homotopy equivalence
\endproclaim

If one assumes Thurston's Geometrization conjectural picture of 3-manifolds,
the answer for (B) is Yes if $k=\infty$ [Ro2]; the answer for (A) is Yes if
the targets are hyperbolic [So], or the domain is non-Haken [RW],
or the targets have finite $\pi_1$ [LMWZ].

Thus it seems natural to study the conjectures of J. Simon
for closed orientable 3-manifolds (Question 1.6 in \S 3).
We find that the positive answer for Question 1.5 are important
for studying the conjectures. This will be addressed in \S 3.

\vskip 0.5 true cm

{\bf Section 2. $\pi_1$-surjective maps between aspherical Seifert manifolds.}

\vskip 0.5 true cm

\proclaim{Theorem 2.1}
Let $M_1$ and $M_2$ be closed orientable
aspherical Seifert fiber spaces with
the same rank and whose base orbifolds are orientable.
Then any $\pi_1$-surjective map
$f: M_1\to M_2$ is of non-zero degree.
\endproclaim

To prove Theorem 2.1, we will make use of
[Ro1], in particular we refer the reader to \cite{Ro1}
for the definition of a {\it vertical pinch} and a {\it squeeze}, and
{\it squeeze torus}.
Call a squeeze is {\it vertical}, if in the squeeze torus, the squeezing
circle meets the regular fiber exactly one point.

Also remember that any orientable Seifert manifold $M$ with orientable
base orbifold of genus $g$ and with $n$ singular fibers
 has a unique
normal form $(g; b; \alpha_1,\beta_1; ... ; \alpha_n,\beta_n)$,
where $0\le \beta_i \le \alpha_i$, $i=1,...,n$.
The orbifold $O_1$ of $M_1$ will denoted by $(g;\alpha_1,...,\alpha_n)$.
$g$ is often omitted if $g=0$.

We first need two lemmas.

\proclaim{Lemma 2.2 [Ro1, Lemma 3.5]} Let $f: M\to N$ be a
 map between aspherical Seifert manifolds
and $1\ne f_*(h)\subset h'$,
where $M$ is closed and $\partial N\ne\emptyset$,
$h$ and $h'$ are regular fibers of $M$ and $N$ respectively.
Then either $f$ admits a vertical squeeze, or $f$ can be homotoped so that
the image of $f$ lies in a fiber of $N$.
\endproclaim

\proclaim{Lemma 2.3}
Suppose $f: F\to O$ is an orbifold branch covering, where $F$ is a surface of
genus $g$, and $O$ is a orbifold,
both are orientable and have non-positive Euler characteristic, then
$rank(\pi_1(F)) \ge rank(\pi_1(O))-1$ if $f$ is a double branched cover
over 2-sphere
and
$rank(\pi_1(F))\ge rank(\pi_1(O))$ otherwise.
\endproclaim

\demo{Proof}
The proof is based on the results about the ranks of
Fuchsian groups [ZVC, Theorem 4.16.1] and the
Riemann-Hurwitz formula.

Suppose $O$ has $k$ singular points of index $v_i$, $i=1,...,k$,
with the underlying space of genus $g'$ and the degree of $f$ is $n$.
Then we have

$$2-2g=n(2-2g'-\Sigma_{i=1}^k(1-\frac 1{v_i}))$$

For the case $g=1$, the verification is direct, so we assume below that $g>1$.

If $n=2$ then all $v_i=2$, $k=2m$ and we have
$2-2g=2(2-2g'-m)$, i.e., $g=2g'+m-1$.
Now $rank(\pi_1(F))=2g=4g'+2m-2$ and the $rank(\pi_1(O))$ is at most $2g'+2m-1$
if $g'>0$ and is $2m-1$ if $g'=0$ by [ZVC, Theorem 4.16.1].
In any case the lemma follows.

If $n\ge 3$, then
$$2-2g\le 3(2-2g'-\Sigma_{i=1}^k(1-\frac 1{v_i}))\le
3(2-2g'-\frac k2) $$
i.e., $2g\ge 6g'-4 + \frac {3k}2$. If $g'>0$, $2g\ge 2g'+\frac {3k}2$.
But the rank of $\pi_1(O)$ is at most $2g'+k-1$.
If $g'=0$, then we have $2g\ge -4+\frac{3k}2 $ if $k$ is even
and $g\ge -4+\frac {3k}2+\frac 12$ if $k$ is odd. The rank of $\pi_1(O)$
is at most $k-1$. It follows that if $k\le 5$, then $2g\ge k-1$. If $k\le 4$
we still have $2g\ge k-1$ since we assume that $g>1$.
\qed

\demo{Proof of Theorem 2.1}
Suppose $f$ is of zero degree.
For clarity, the proof is divided up three steps.

{\bf {Step (1)}}
We prove the following:

Claim : $f(h)$ is homotopically non-trivial,
where $h$ is the regular fiber of $M_1$.

\leftline{\bf Proof of Claim:}

Let $M_1=(g;b; a_1,b_1;...;a_k,b_k)$ and $G_1=\pi_1(M_1)/<h>$,
where $<h>$ is the cyclic group generated by regular fiber of $M_1$.
Let $r=rank(\pi_1(M_1))=rank(\pi_2(M_2))$.

By [BZ, Theorem 1.1] and [ZVC, Theorem 4.16.1], one of the following cases
holds.

(1) $rank(\pi_1(M_1))> rank(G_1)$,

(2) $rank(\pi_1(M_1))=rank(G_1)$ there is  a set of generators of $G_1$ which
realizes the rank and contains at least one  torsion element,

(3) $rank(\pi_1(M_1))=rank(G_1)=rank(G_1/T)$, where $T$ is the normal
subgroup normally generated by the torsion elements and $G_1/T$ is a surface group.

If $f(h)$ is homotopically trivial,
then $f_*: \pi_1(M_1) \to \pi_1(M_2)$ induces
an epimorphism $\phi: G_1\to \pi_1(M_2)$.

In Case (1) the Claim is clearly true.

In Case (2) the Claim is also true  since $\pi_1(M_2)$ is torsion free.

In Case (3), $f_*$ induces an epimorphism $\phi' : G_1/T \to \pi_1(M_2)$.
Let $f' : F\to M_2$ be the map which realizes $\phi'$.
Since $\phi'$ is not injective,  (otherwise $\phi'$ would be an isomorphism
and $\pi_1(M_2)$ would be surface group),
by the simple loop theorem for maps from a surface to a Seifert manifold [H],
there are essential simple  loops in the kernel of $\phi'$.
Assume first there is an essential non-separating simple loop, which we
$\alpha$, in the kernel.
Then the map $f'$ induces a map $f'':F' \to M_2$,
where $F'$ is a complex obtained by squeezing $F$ along   $\alpha$.
It is easy to see that the rank of $\pi_1(F')$ is $r-1$. We reach
a contradiction. If all essential simple loops in kernel of $\phi'$ are
separating, let $\alpha$ be a maximal family of non-parallel
separating essential simple closed curves in kernel $\phi'$.
Again $f'$ can factors through
$f'':F' \to M_2$,
where $F'$ is a complex of obtained by squeezing $F$ along   $\alpha$,
which is union of closed surfaces connected by arcs.
Let $S$ be a surface in $F'$.
Due to the maximality of $\alpha$,
the restriction $f''|_{S}$  $\pi_1$-injective,
which must be either horizontal or vertical by [H].
If $f''|_{S}$ is horizontal, than $p_2\circ f''|: S\to O(M_2)$
is an orbifold branched covering, where $p_2 :M_2 \to O(M_2)$ is the fiber map.
But the rank of $\pi_1(S)$ is at most $r-2$.
This is also ruled out by Lemma 2.3.
If $f''|_{S}$ is vertical for each surface $S$ of $F'$,
then $F'$ contains at most $g$ such surfaces and and each of them is a torus.
Clearly the rank of $f''_*\pi_1 (F')$ is at most $g+1$, which is
at most $r-1$ (since $g>1$ and
$r\ge 2g$). Again we reach a contradiction.

\leftline{\bf \hbox{Step (2)}}
We will factor $f: M_1\to X\to M_2$, where the 2-dimensional
complex $X$ is a quotient of $M$ with rank $r_X$.

Since $f(h)$ is homotopically non-trivial, a standard argument in
3-manifold topology shows that $f:M_1\to M_2$ can be deformed to be a
fiber preserving map (see [J] for example).
Since a vertical pinch reduces the rank of
$\pi_1$, $f:M_1\to M_2$ admits no vertical pinch.  Suppose the mapping
degree is zero. We can further deform the map so that the image
$f(M_1)$ misses a regular fiber $h'$ of $M_2$.  To see this, 
$f: M_1\to M_2$ is fiber preserving. We can further deform $f$
so that for each singular fiber of $M_2$, its preimage consists
of finitely many fibers of $M_1$.
Now remove all singular fibers of $M_1$ and their f-images,
and remove all singular fibers of $M_2$ and their f-preimages.
The restriction of $f$ gives a proper map $f': M'_1\to M'_2$,
which is fiber preserving map between circle bundles. 
Since $f$ is assumed to be degree zero, $f'$ is of zero degree.
Since $f(h)$ is non-trivial, the induced proper map $\bar f': F'_1\to F'_2$
between base surfaces must be degree zero. Hence by
$\bar f'$ can be deformed so that its image misses
a point of $F'_2$. This deformation can be lifted to the bundle map $f'$
whose image then misses a circle fiber in $M'_2$.
With this we reach the situation claimed above.

Now remove an open fibered
neighborhood of $h'$, and denote the resulting manifold by $N$. Then
we have a fiber preserving map $f: M_1\to N$, where $\partial N\ne
\emptyset$.

According Lemma 2.2
either $f: M_1\to N $ admits a fiber squeeze along an incompressible
vertical torus, or $f(M_1)\subset $ a fiber of $N$.
Using this we can reformulate the above so that either $f: M_1\to M_2$
admits a vertical squeeze along an incompressible
vertical torus, or $f(M_1)\subset $ a fiber of $M_2$.

Since $f$ is $\pi_1$-surjective, and $M_2$ is an closed aspherical
Seifert fiber space, the situation that $f(M_1)\subset $ a fiber of $M_2$
cannot happen. Let $\Cal T$ be a maximal family of
disjoint non-parallel incompressible tori along
which $f$ admits vertical squeeze.
Let $X=\Cal Q\cup \Cal A$ be the space obtained after the squeezing,
where $\Cal Q$ is a union of
Seifert fiber spaces with the induced Seifert
fibration, $\Cal A$ is a union of annuli and $\partial \Cal A$
is a union of regular fibers  of $\Cal Q$.
Then $f$ induces a $\pi_1$-surjective map $X \to M_2$,
which we continue to denote by $f$.

Now all components of $\Cal Q$ are Seifert fibered spaces with the induced
Seifert fibrations, so 
we may assume that $Q_1,...., Q_{k_1}$ are Seifert manifolds
of $\Cal Q$ which are not the trivial circle bundle over $S^2$ and
$Q_{k_1+1},....,Q_{k_1+k_2}$
are  trivial circle bundle over $S^2$.
Clearly each $Q_j$, $j>k_1$, is $S^2\times S^1$.

For each $j>k_1$, there is an annulus $A$ in $\Cal A$,
with two components $C_1$ and $C_2$ of $\partial A$
such that $C_1$ belongs to $Q_j$,
Since $C_1$ is a regular fiber of $Q_j=S^2\times S^1$,
 $Q_j\cup A$ has $C_2$ as a  retractor,
 and hence we can send $Q_j\cup A$ to $C_2$ by this retract,
then extend the map to whole $\Cal Q\cup \Cal A$.
After $k_2$ such operations, we get a quotient
space $X_1=\Cal Q_1\cup\Cal A_1$
where $\Cal Q_1=\{Q_1,...,Q_{k_1}\}$,
and $f$ induces a $\pi_1$-surjective map $X_1\to M_2$,

By the maximality of $\Cal T$,
each $Q_i$ contains no squeeze torus for $f|_{Q_i}$, so
we have that $f(Q_i)\subset$
a fiber of $M_2$, and consequently we have the following

\leftline{\bf Fact~(*)}: each $Q_i$ has base orbifold $S^2$
and has no more than 3 singular fibers
(otherwise there will be a squeeze torus).

So $f :X_1\to M_2$  induces a
$\pi_1$-surjective map $X={\Cal S} \cup {\Cal A_1}\to M_2$,
$\Cal S$ is a union of $k_1$ circle.

\leftline{\bf \hbox{Step (3)}}
We will show that $r_X < r$ and reach a contradiction.

If $g>0$, there is a non-separating
squeeze torus for $f$ and clearly $r_X<r$.

Below we assume that $g=0$. Then every squeeze torus is a separating torus.

Say that $M_1$ is of type I, if
$M_1$ has normal form $(0; b; 2,1;....;2,1;2\lambda+1,b_k)$,
where $k\ge 4$ is even, and $\lambda>1$,
otherwise call $M_1$ is of type II.
By [Theorem 1.1 BZ], the $r=k-2$ if $M_1$ is of type I and
$r=k-1$ if $M_1$ is of type II.

Each component $Q_i$ of $\Cal Q$ must have infinite
fundamental group, otherwise $f(h)$ is an element of finite order,
which must be trivial in
$\pi_1(M_2)$, and this is forbidden by Step (1).
In particular, each $Q_i$ contains at least 2 singular fibers,
$i=1,...,k_1$, and $Q_i$ contains exactly
two singular fibers only if $Q_i=(0;0;a,b;a,-b)$.
Moreover if $M_1$ is of type I, then at least one $Q_i$ contains
4 singular fibers (since $\lambda>1$ and both
$(0;b;2,1;2\lambda+1, b_2)$ and
$(0;b;2,1;2,1;2\lambda+1, b_3)$ have finite fundamental groups), which is not possible
by Fact*.

If $M_1$ is of type II, then $r=k-1$ and $k\ge 3$, but
 $$r_X \le k_1 \le \frac k2< k-1=r,$$
where the first $\le$ is due to the fact that every squeeze torus is
separating and the second $\le$ is due to every $Q_i$ contains
at least two singular fibers.
\qed

{\bf Remark 2.4.}
In Theorem 2.1, the condition "$f$ is $\pi_1$-surjective"
can be replaced by "$f$ is $\pi_1$-finite-index",
and the condition "orbifolds are orientable" can be removed.
For details see [Hu], where the proof is parallel to the proof
above, but involves more complicated case by case argument.

{\bf 3. On  the conjectures of J. Simon on 3-manifold groups.}

\vskip 0.5 true cm

In this section we study the follwing questions.

\proclaim{Question 1.6} Let $M_i$ be closed
orientable aspherical 3-manifolds.
Suppose there is an epimorphism $\phi : \pi_1(M_1)\to \pi_1(M_0)$.

(A) Is rank $\pi_1(M_1)>$ rank $\pi_1(M_0)$?

(B) Is Heegaard genus of $M_1\ge$ Heegaard genus $M_0$?

Moreover given $M_0$.

(C1) Is there a number $N_M$ such that any sequence
of epimorphisms
 $\pi_1(M_0)\to \pi_1(M_1)\to .... \to \pi_1(M_n)$ with $n\ge N_M$ contains
an isomorphism?

(C2)
Does any infinite sequence of epimorphisms
 $\pi_1(M_0)\to \pi_1(M_1) $
 $\to .... \to \pi_1(M_n)\to ...$ contain
an isomorphism?

(D) Are there only finitely many  $M_i$ with the same first
Betti number, or the same $\pi_1$-rank, as that of $M_0$, for which
there is an epimorphism $\pi_1(M_0)\to \pi_1(M_i)$?
\endproclaim

We remark that a positive answer for (B) of Question 1.6 implies a
positive solution to the Poincare Conjecture. From Example 1.4 of The 
Introduction
the answer to (D) is negative if we remove the condition on first Betti
number or $\pi_1$-rank on (D) of Question 1.6.

We describe first some examples of non-trivial
$\pi_1$-surjective maps
between two 3-manifolds of the same rank,
which give a negative answers of the (A) of Question 1.6.  Clearly
those examples are all of non-zero degrees.

{\bf {Example 3.1.}}
Let $M$ be a Seifert manifold of normal form $(0;0;6, b_1; 5, b_2; 7, b_3)$.
Let ${\Bbb Z}_2$ be a cyclic group acting on $M$ such that
it induces the identity on the base space and standard rotation
on each regular fiber. Then one verifies that
$M/{\Bbb Z}_2$ is a Seifert manifold with normal form
$(0;0;3, b_1; 5, 2b_2; 7, 2b_3)$. Now
$$\pi_1(M)=<s_1, s_2,s_3, h~|~[s_j,h], s_1^6h^{b_1},s_2^5h^{b_2},s_3^7h^{b_3},
s_1s_2s_3>$$
and
$$\pi_1(M/{\Bbb Z}_2)=<t_1, t_2,t_3, h'~ |~ [t_j,h'],
t_1^3{h'}^{b_1},t_2^5{h'}^{2b_2}, t_3^7{h'}^{2b_3}, t_1t_2t_3>$$ The
quotient map $p:M\to M/{\Bbb Z}_2$ is a branched covering of degree 2 and
$p_*$ sends $s_j\mapsto t_j$ and $h\mapsto {h'}^2$. Since $(2, b_1)=1$,
$p_*$ is surjective. By [BZ] these manifolds have rank 2.
$\square$

{\bf Examples 3.2.}
We now give some examples of $\pi_1$-surjective non-zero
degree maps between hyperbolic manifolds of the same $\pi_1$ ranks.

Let $M$ be a closed orientable 3-manifold and $k \subset M$ 
be any hyperbolic fibered knot.
Suppose the fiber $F$ has genus $g$.
Let $M_n$ be the $n$-fold cyclic branched cover of $M$ over the knot $k$.
Then the rank of $\pi_1(M_n)$ is bounded by $2g+1$ for all $n$
and $M_n$ is hyperbolic when $n$ is large.
If $k|n$, then $M_n\to M_k$ is a branched cover,
which is $\pi_1$-surjective.
So there are must be infinitely many $\pi_1$-surjective branched covering
$M_n \to M_k$ between hyperbolic 3-manifolds of the same ranks.

A well studied case is when $M_n$ is the n-fold
cyclic branched cover of the figure eight knot.
 Then  for $n\geq 3$ the fundamental
 groups are all 2-generator---in fact they are 
the Fibonacci groups $F(2,2n)$, which are all hyperbolic if $n\ge 4$.
By abelianizing $F(2,2n)$ we see that all $M_n$ have first Betti
number zero (see [MR] for example).

The next example gives the negative answer of (C1) of Question 1.6.

{\bf Example 3.3}

(1) Let $M(n,k)=(0;0; 2^k3,b_1; 5, 2^{n-k}b_2; 7, 2^{n-k}b_3)$.
Similar to Example 3.1, we have sequence of degree 2 branched covering
$M(n,n)\to ...\to M(n,1)\to M(n,0)$
of length $n+1$,
which induces a sequence of epimorphisms of groups
$\pi_1(M(n,n))\to ... \to \pi_1(M(n,1)))\to \pi_1(M(n,0))$
of rank 2. Let $M$ be $\Sigma_2\times S^1$.
Clearly $\pi_1(M)$ surjects onto $Z* Z$, then we
have the sequence of epimorphisms
$$\pi_1(M)\to\pi_1(M(n,n))\to ... \to \pi_1(M(n,1))\to \pi_1(M(n,0))$$
of length $n+2$, where $n$ can be arbitraily large.

Moreover if we choose $b_1, b_2, b_3$ such that the Euler number
of $M(n,n)$ is non-zero. Since each $M(N,k)$ has infinite $\pi_1$ and
is the image of $M(n,n)$ under non-zero degree map,
the Euler number of $M(n,k)$ is non-zero [Theorem 2, W]. It follows
$M(n,k)$ has neither horizontal or vertical incompressible surface,
and therefore all $M(n,k)$ are non-Haken [J].

(2) Let $M_n$ be the $n$-fold cyclic branched covering
of $S^3$ over figure eight knot as in the end of Example 3.2.
Then we have sequence of branched coverings of
hyperbolic rational homology spheres
$M_{4k}\to ...\to M_8\to M_4$
of length $l$ which induces a sequence of epimorphisms of groups
$\pi_1(M_{4k})\to ... \to \pi_1(M_8)\to \pi_1(M_4)$ with rank 2.
Let $M$ be a hyperbolic 3-manifold with $\pi_1(M)$ surjecting 
${\Bbb Z} * {\Bbb Z}$ (as in Example 1.4). Then we
have the sequence of epimorphisms
$$\pi_1(M)\to\pi_1(M_{4k})\to ... \to \pi_1(M_8)\to \pi_1(M_4)$$
of length $l+1$, $l$ can be arbitraily large. \qed

The next result gives a partial positive answer of (C2) of Question 1.6.

\proclaim{Theorem 3.4}
Given $M_0$, and a sequence $M_i$ of closed
orientable aspherical Seifert manifolds with epimorphisms
$\pi_1(M_0)\to \pi_1(M_1)\to .... \to \pi_1(M_n)$ $\to ...$, this
sequence contains an isomorphism.
\endproclaim

\demo{Proof}
By passing an infinite subsequence, we may assume that
 all groups in the sequence have the same rank (each epimorphism in
 the subsequence is the composition of epimorphisms involved).
Then each epimorphism $\phi_i: \pi_1(M_i)\to \pi_1(M_{i+1})$
in the sequence can
be realized by a map $f_i: M_i\to M_{i+1}$ of non-zero degree
by Theorem 2.1.
Moreover the Seifert fibrations of the $M_i$'s can be arranged so that
each  $f_i$ is a fiber preserving.
Let $O_i$ be the orbifold of $M_i$, then $\chi(O_i)\le 0$ and we
have the induced sequence of epimorphisms

$$\pi_1(O_0)\to \pi_1(O_1)\to .... \to \pi_1(O_n)\to...$$
of Fuchsian groups.
We therefore have a decreasing sequence

$$\chi(M_0)\le \chi(M_1)\le .... \chi \pi_1(M_n)\le....$$
The $\{-\chi(O)\}$ form a well-order subset of reals,
where $O$ runs over compact orbifolds,
$\chi (O_k)=\chi(O_{k+1})$ for $k$ larger than a given $N$
([Ro2, Lemmas 2.5 and 2.6] for details).
Since there are at most finitely many orbifolds $O$ with given
$\chi$, by passing an infinite sequence, we may assume that all
$O_i$ are the same.

Let  $O_i=(g;\alpha_{1,}, ...,\alpha_{n})$.
Then $M_i=(g; b_i; \alpha_{1,}, \beta_{1,i};...;\alpha_{n,},\beta_{n,i})$.

Since $0<\beta_{l,i}<\alpha_{l}$ for $l=1,...n$,
by passing a further subsequence, we may assume that
$\beta_{l,i}=\beta_l$, and finally we get
$M_i=(g; b_i; \alpha_1, \beta_1;...;\alpha_n,\beta_n)$.
Moreover we may assume that all $b_i\ne 0$.
Note that by [p. 680 of LWZ], all $M_i$
have the same first Betti number and the torsion part of $H_1(M_i,{\Bbb Z})$
is unbounded if $b_i$ unbounded.
Since  epimorphisms on $\pi_1$ induce epimorphisms
on first homology groups, it follows that $b_i$'s are bounded.
Now we have $b_i=b_j$ for some $i,j$, then $M_i=M_j$ and
by the hopfian property of Seifert
manifold groups, the epimorphism $\pi_1(M_i)\to \pi_1(M_j)$ is an
isomorphism. Then in the sequence above
there must be an isomorphism. Theorem 3.4 follows.\qed

We have seen that Theorem 2.1 plays important roles for the proof Theorem
3.4. If the  answer of
Question 1.5 is also YES for hyperbolic 3-manifolds,
this will lead to a positive answer for 
(C2) and (D) for hyperbolic 3-manifolds.

\proclaim{Proposition 3.5}
Suppose Question 1.5 has a positive answer for hyperbolic 3-manifolds.
Then for a given closed orientable hyperbolic 3-manifold $M_0$,

(1) any infinite sequence of epimorphisms
 $\pi_1(M_0)\to \pi_1(M_1)\to .... \to \pi_1(M_n)\to ...$ contains
an isomorphism, where all $M_i$ are closed orientable  hyperbolic 3-manifolds.

(2) there are only finitely many closed orientable hyperbolic 3-manifolds
$M_i$ with the same $\pi_1$-rank as that of $M_0$, for which
there is an epimorphism $\pi_1(M_0)\to \pi_1(M_i)$.
\endproclaim

\demo{Proof} (1) By passing an infinite subsequence we may assume all
$\pi_1(M_i)$ have the same rank. Since we assumed that Question 1.5 
has a positive answer for hyperbolic
3-manifolds, this sequence is realized by
a sequence of non-zero degree maps
$$M_0\to M_1 \to ... \to M_n\to...$$

The rest of the proof is now standard. Since all maps $f_i: M_i\to
M_{i+1}$ in the sequence are of non-zero degree, by Gromov's Theorem
[Chapter 6, Th], $v(M_i)\ge v(M_{i+1})$, where $v(M_i)$ is the
hyperbolic volume of $M_i$.  By Thurston-J\o genson's Theorem [Chapter
6, Th], $v(M_k)$ must be a constant when $k$ is larger than a given
integer $N$. Then by Gromov-Thurston's Theorem [Chapter 6, Th], $f_k$
is homotopic to a homeomorphism, $k>N$, so $f_{k*}$ is an isomorphism

For (2) since we again assume that Question of 1.5 has a
positive answer for hyperbolic 3-manifolds, each $\phi_i: \pi_1(M_0)\to
\pi_i(M_i)$ can be realized by a map of non-zero degree. By Soma's theorem
[So], there are only finitely many such $M_i$.
\qed

We also note the following partial positive answer of (D) of Question 1.6
follows easily from the methods of [RW].

\proclaim{Theorem 3.6} Suppose $M$ is a non-Haken hyperbolic 3-manifold.
Then there are are only finitely many  closed orientable hyperbolic
3-manifolds $M_i$  for which
there is an epimorphism $\pi_1(M)\to \pi_1(M_i)$.
\qed\endproclaim

\vskip 0.5 true cm

{\bf Section 4. A minimal Haken manifold}

\vskip 0.5 true cm

Let $E$ be the complement of trefoil knot with
$m$ the meridian and $l$ the longitude.
$E$ has a unique Seifert fibration with two singular fiber of indices
2 and 3, over the disc. Via this Seifert structure, we have a presentation
$$\pi_1(E)=<a, b, c, t~|~a^2t, b^3t, abc>$$
where $t$ is the regular Seifert fiber.
Let $E_1$ and $E_2$ be homeomorphic to $E$ with meridian and longitudes
$(m_i,l_i)$, $i = 1,2$.
Now glue $E_1$ to $E_2$ via a homeomorphism $h: \partial E_1\to \partial E_2$
such that $h(l_1)=m_2$ and $h(m_1)=l_2^{-1}$.
Let $M$ denote the resulting manifold, which is a closed graph manifold.
The main theorem of this section is:
\proclaim{Theorem 4.1}
$M$ is a minimal closed Haken 3-manifold
among all 3-manifolds satisfying Thurston's geometric conjecture.
\endproclaim
We begin the proof by collecting some elementary facts.
\proclaim{Lemma 4.2}
(1) For any representation $\phi : \pi_1(E)\to SL(2,{\Bbb C})$,
if $\phi (t)\ne 1$, then the image $\phi(\pi_1(E))$ is a cyclic group
$<\lambda>$. Moreover, we must have $\phi(a)= \lambda^{-2}$,
 $\phi(b)= \lambda^{-3}$, $\phi(c)= \lambda^{5}$,  and $\phi(t)=\lambda^{6}$.

 (2) In $\pi_1(T)$, where $T=\partial E$,
  we have $m=tc^{-1}$ and $l=t^{-5}c^{6}$. (Equivalently,
 $t=6m+l$ and $c=5m+l$.) Hence  $h(t_1^{-5}c_1^{6})=t_2c_2^{-1}$ and
 $h(t_1c_1^{-1})=t_2^5c_2^{-6}$.

 (3) $M$ is an integral homology 3-sphere.

 (4) the only 2-sided incompressible surface is the incompressible
 torus $T$, which separates $M$ into $E_1$ and $E_2$.
 \endproclaim

 \demo{Proof}
 The main part of (1) follows from \cite{M, Prop. 3} and the fact that
 $H_1(E,{\Bbb Z})$ is cyclic.
 (2) and (3) and the remaining parts of (1) are just direct calculations.
 Finally to establish (4) we observe the following. Since the trefoil knot
 is 2-bridge $E$ cannot contain a closed embedded essential surface by
 \cite{HT}. If $M$ contained an embedded incompressible surface $\neq T$,
 it would follow from the remark above and the gluing homeomorphism that
 $E$ would have a boundary slope $1/0$. However \cite{Theorem 2.0.3, CGLS} then implies
 the existence of a closed embedded essential surface in $E$.
 \qed

 To show that $M$ is minimal, we assume not and suppose that there is
 a degree one map $f:M\to N$, where $N$ is irreducible,
 $N\ne M$, and  $N\ne S^3$.
 First, since $M$ is a graph manifold, its Gromov norm is zero,
 so $N$ cannot be a hyperbolic 3-manifold by \cite{T, Chapter 6}.
 Moreover it is well-known that $N$ must be an integer homology sphere
 ([Lemma 3.1 RW].
 The proof of Theorem 4.1 will be finished by Lemmas 4.3,
 4.4 and 4.5 below.

 \proclaim{Lemma 4.3}
 $N$ is non-Haken.
 \endproclaim

 \demo{Proof}
 Suppose $N$ is Haken, and let $F\subset N$ be an embedded
 incompressible surface.  We may deform $f$ so that $f^{-1}(F)$ is an
 incompressible surface in $M$.  By (4) of Lemma 4.2
 $f^{-1}(F)$ must consist of parallel copies of $T$.  By standard
 3-manifold topology, we can further deform $f$ so that
 $f^{-1}(F)=T$. It follows that $F$ is a 2-sphere or torus.
 Since $N$ is irreducible, $F$ must be a torus separating $N$ into
 two parts $N_1$ and $N_2$. Furthermore, the map $f$ can be decomposed into two
 proper degree one map $f|: E_i\to N_i$.  However $E_i$ is a minimal
 3-manifold among knot complements in 3-manifolds via proper degree one
 maps [BW]. Thus, each $f|$ is a homeomorphism, and it follows that $f$
 itself is homotopic to a homeomorphism.
 \qed

 \proclaim{Lemma 4.4}
 $N$ is not a Seifert manifold with finite fundamental group (other than
possibly $S^3$).\endproclaim

 \demo{Proof}
 By (3) of Lemma 4.2 if $N$ is a Seifert fibered manifold
 of finite fundamental group and $N\ne S^3$, it must be the Poincar\' e
 Homology 3-sphere $P$. Note $\pi_1(P)$ surjects onto $A_5$, the 
alternating group on 5 letters. In particular,
 (as is well-known) $A_5$ is a subgroup of $PSL(2,{\Bbb C})$---since 
 $SO(3)$ can be
 identified with $PSU(2)$, and the latter is a subgroup
 of $PSL(2,{\Bbb C})$.
 To prove the lemma, it suffices to prove that the image group of any
 representation of $\phi : \pi_1(M) \mapsto PSL(2,{\Bbb C})$ cannot be $A_5$.

 Case (1) If $\phi(t_1)\ne 1$ and $\phi(t_2)\ne 1$, by (1) of Lemma 4.2,
 the whole image $\phi(\pi_1(M))$ must be a cyclic group (actually trivial).

 Case (2) Without loss of generality, we may assume that $\phi(t_1)=1$
 and $\phi(t_2)\ne 1$. By (1) and (2) of Lemma 4.2,
 $\phi: \pi_1(M)\to PSL(2,{\Bbb C})$ factors as $\nu: \pi_1(M)\to G$
 and $\mu: G\to PSL(2,{\Bbb C})$ where $G$ is generated by two groups
 described in (a) and (b) below:
 
 (a) $\nu(\pi_1(E_1))=<a_1, b_1, c_1| a_1^2, b_1^3, a_1b_1c_1>$,
 (b)  A cyclic group $<\lambda_2>$ such that
  $\nu(c_2)= \lambda_2^5$,  $\nu(t_2)= \lambda_2^6$.

 Since  $h(t_1c_1^{-1})=t_2^5c_2^{-6}$, we have $\nu(h(c_1^{-1}))
 =\nu (h(t_1c_1^{-1})) = \nu t_2^5c_2^{-6}) = 1$.
 It follows that
 $$G=<a_1, b_1, c_1| a_1^2, b_1^3, a_1b_1c_1, c_1> =
 <a_1, b_1| a_1^2, b_1^3, a_1b_1>,$$
 which is the trivial group.

 Case (3)  $\phi(t_1)= 1$ and $\phi(t_2)= 1$.
 In this case $\phi: \pi_1(M)\to PSL(2,{\Bbb C})$ factors through a group $G$
 via a map $\nu: \pi_1(M)\to G$, with
 $\nu(\pi_1(E_i)$ is the quotient of
 $G_i=<a_i, b_i, c_i| a_i^2, b_i^3, a_ib_ic_i>$, $i=1,2$.
 Moreover by (2) of Lemma 4.2 we have that in the quotient
 $c_1=c_2^6$ and $c_2=c_1^{-6}$. Immediately we have
 that $c_1^{37}=1$ and  $c_2^{37}=1$
and finally
$$G=<a_i, b_i, c_i, i=1,2| a_i^2, b_i^3, a_ib_ic_i, c_i^{37},c_1=c_2^6,
i=1,2>$$
Suppose there is a homomorphism $\mu :G\to A_5$. Since the order of $c_i$
is 37, and $A_5$ has order $60$, under the homomorphism $\mu$
the images of $c_1$ and $c_2$ must be trivial.
It follows that  $\mu :G\to A_5$ can factor through
the group $G'$,
$$G'=<a_1, b_1,| a_1^2, b_1^3, a_1b_1>*<a_2, b_2,| a_2^2, b_2^3, a_2b_2>,$$
but as above this is trivial.
\qed

\proclaim{Lemma 4.5}
$N$ is not a Seifert manifold with infinite $\pi_1$.
\endproclaim
\noindent The proof of Lemma 4.5 requires
a sequence of additional lemmas.
Suppose below $N$ is a Seifert manifold of infinite $\pi_1$. By Lemma 4.3,
we may assume that $N$ is non-Haken. Then $N$ must be a Seifert manifold
with three singular fibers over $S^2$.

We begin by establishing:

\proclaim{Lemma 4.5.1}
(1) Suppose  $\Delta\subset Iso_+H^2$ is a triangle group
and $\phi: \pi_1(2,3,l)\to \Delta$ is of finite index.
Then  the image of $\phi$ is a hyperbolic triangle group
isomorphic to $\pi_1(2,3,k)$, where $k|l$.

(2) Suppose a Serfert manifolds $N$ is an integer homology sphere
with infinite $\pi_1$ and orbifold $O=(a_1,a_2,a_3)$. Then
$gcd(a_i, a_j)=1$ for $i, j=1,2,3$, and $O$ is a hyperbolic
orbifold.
\endproclaim
\demo{Proof}
(1) Let $x', y'$ be the order 2 and order 3 elements which generates
$\pi_1(2,3,l)$ such that $x'y'$ is of order $l$. Use $x$ and $y$ to denote
their images in $Iso_+H^2$, then $x$ and $y$ generate the image of $\phi$.
Since the image of $\phi$ is of finite index in $\Delta$, it must be
co-compact and of rank 2. By well-know fact then the image is a triangle group
with $x^2=y^3=(xy)^k=1$, where $k|l$.

(2) follows from [p. 680 (d) LWZ].
\qed
\proclaim{Lemma 4.5.2}
There is a simple closed curve in the kernel of $f| : T\to N$.
\endproclaim
\demo{Proof}
Since $\pi_1(N)$ is torsion free
and $T$ is a torus,
to prove the lemma, we need only that
the kernel of $f| : T\to N$.

Suppose $f(t_1)\ne 1$, otherwise the claim is proved.
Note that all elements in $f(\pi_1(E_1))$ commute with $f(t_1)$.
If $f(t_1)$ is not the fiber $t$ of $N$,  then either
$$f(\pi_1(E_1)) = f(t_1)\, \hbox{or}\, f(\pi_1(E_1)) = <f(t_1), f(c_1)> =
{\Bbb Z}\oplus {\Bbb Z}.$$
The second case is not possible since $H_1(E_1;{\Bbb Z})={\Bbb Z}$.  
In the first case
we deduce that $\ker (f|_T)_*$ is nontrivial.  Similarly if $f(t_2)$
is not the fiber $t$ of $N$, then $\ker (f|_T)_*$ is non-trivial.
If $f(t_1)=t=f(t_2)$. Since $t_1$ and $t_2$ do not coincide up to isotopy,
still we have $ker (f|_T)_*$ is non-trivial.
\qed

Let $C$ be the simple closed curve provided by Lemma 4.5.2.
Suppose $C=pm_1 + ql_1$ on $\partial E_1$, then
$C=-qm_2+pl_2$. By (1) of Lemma 4.2 we have
$pm+gl=(p-5q)t+(-p+6q)c$ and $-qm+pl=(-q-5p)t+(q+6p)c$.
So the degree 1 map $f$ factors through $f: M\to N_1\cup_{S^1} N_2
\to N$ where $N_1$ and $E_2'$ are Seifert manifolds whose normal
forms are given by $(2,1; 3,1; -p+6q, p-5q)$ and $(2,1; 3,1; q+6p,-5p-
q)$ respectively, and the two cores of the surgery solid tori are identified.
If $f|_*(\pi_1(N_1))\ne \pi_1(N)$ and $f|_*(\pi_1(N_2))\ne \pi_1(N)$
then $\pi_1(N) $ can be presented as a non-trivial free product with
amalgamation by the classical result (see [CGLS] for example).
It follows that $N$ will is Haken contrary to Lemma 4.3.
Thus without loss, we assume that $f|_*(\pi_1(N_1))= \pi_1(N)$.

\proclaim{Lemma 4.5.3}
$f|_{N_2}$ is of degree non-zero.
\endproclaim

\demo{Proof}
Let $\tilde E$ be the covering of $N$ corresponding to $f|_*(\pi_1(N_2))$.
Then $f: N_2\to N$ lifts to $\tilde f: N_2\to \tilde E$,
which is $\pi_1$-surjective.
If $f|_*(\pi_1(N_2))\subset \pi_1(N)$ is of finite index,
then $\tilde E$ is a closed Seifert manifold.
Since both $\pi_1(N_1)$ and $\pi_1(N)$ are rank 2,
$\pi_1(\tilde E))$ must be also rank 2.
Then $\tilde f$ is of non-zero degree by
Theorem 2.1. Hence $f|_{N_2}$ is non-zero degree.

Below we show  $f|_*(\pi_1(N_2))\subset \pi_1(N)$
must be of finite index. Otherwise $\tilde E$
is a non-compact, aspherical Seifert manifold, which is known
that either the rank of $H_1(\tilde E)$ is positive
or  $\pi_1(\tilde E)$ is trivial. Since $f|_*(\pi_1(N_2))$ is not trivial
and $N_2$ is a rational homology sphere,
all of the above cases are ruled out.
So $f|_*(\pi_1(N_2)$ must be of finite index in $\pi_1(N)$.
\qed

Since $N_1$ and $N_2$ are in symmetry position,
we have both $f|{N_1}$ and
$f|{N_2}$ are of non-zero degree.

By  Lemma 4.5.3,
we may assume that $f|N_i$ is fiber preserving. Then $f|N_i$ induces an
homoporphism $\phi_i :\pi_1(O_i)\to \pi_1(O)$,
in particular $\phi_1$ is surjective and $\phi_2$
is finite index, where $O_1=(2,3,6q-p)$, $O_2=(2,3,6q+p)$
and $O=(a_1, a_2, a_3)$ are orbifolds of $N_1$, $N_2$ and $N$ respectively.
Since $N$ is an integer homology sphere of infinite $\pi_1$,
it follows that $\pi_1(O)$ is isomorphic to a hyperbolic triangle group.
Since $\phi_1: G_1\to G$ is surjective, it follows that $O=(2,3,k)$,
where  $k|6q-p$ by Lemma 4.5.1 (1).
Since $\phi_2$ is of finite index,
the image of $\phi_2$ is a hyperbolic triangle groups $\pi_1(2,3,k')$
with $k'|6q+p$ by Lemma 4.5.1 (1), moreover $k'|k$.
It is easy to see that $k'$ is a dvisor of both $12q$ and $2p$.
Since $p$ and $q$ are coprime, the great common divisor of $12q$ and $2p$
is 12. So $k'$ is either  2,  or 3, or 4, or 6, or 12.
Then $N$ can not be an integer homology sphere by 4.5.1 (2).
 \qed
\vskip 0.5 true cm

{\bf REFERENCES}

\vskip 0.5 true cm

 [BW] M. Boileau and S. C. Wang {\it Non-zero degree maps and surface
 bundles over $S^1$},  J. Diff. Geom. {\bf 43} (1996), 789--908.

 [BZ] M. Boileau and H. Zieschang, {\it Heegaard genus of
 closed orientable Seifert 3-manifolds}, Invent. Math. {\bf 76} (1984)
 455--468.

 [CGLS] M. Culler, C. McA. Gordon, J. Luecke and P. B. Shalen
 {\it Dehn surgery on knots}, Annals of Math. {\bf 125} (1987) 237--300.

 [CR] P. J. Callahan and A. W. Reid, {\it Hyperbolic
structures on knot complements}, Chaos, Solitons and Fractals, {\bf 9}
(1998), 705--738.

 [H] J. Hass, Minimal surfaces in manifolds with $S^1$ actions and the
 simple loop conjecture for Seifert conjecture for Seifert fiber spaces.
 Proc. AMS, 99, 383-388 (1987)

 [HT] A. Hatcher and W. Thurston. {\it Incompressible surfaces in
 2-bridge knot complements.} Invent. Math. {\bf 79} (1985) 225--246.

 [Hu] C. Huang, Master Thesis, Peking Univ. (1998)

 [J] W. Jaco, Topology of three manifolds, Regional Conference Series
 in Mathematics. 43, AMS Math. Soc. Providence, RI 1980.

 [K] R. Kirby, Problems in low dimensional topology, Geometric Topology,
 Edited by H. Kazez, AMS/IP 1997.

 [LWZ] C. Hayat-Legrand, S. C. Wang and H. Zieschang, {\it Minimal
 Seifert manifolds},  Math. Annalen. {\bf 308} (1997), 673--700.

 [MR] C. Maclachlan and A.W. Reid, Generalised Fibonacci Manifolds,
 Transformation groups, {\bf 2} (1997), 165--182.

 [M] K. Motegi, {\it Haken manifolds and representations of
 their fundamental groups in $SL(2,{\Bbb C})$}, Topology and its
 Applications {\bf 29} (1988), 207--212.

 [RW] A.W. Reid and S. C. Wang,
 {\it Non-Haken 3\--manifolds are not large with respect to mappings of
 non-zero degree}, to appear in Comm. Analysis $\&$ Geom.

 [Ro1] Y. Rong, {\it Maps between Seifert fibered spaces of infinite
 $\pi_1$}, Pacific J. Math. {\bf 160} (1993), 143--154.

 [Ro2] Y. Rong, Degree one map between geometric 3-manifolds, Trans.  A. M. S.
 {\bf 332} (1992), 411--436.

 [Sc] G. P. Scott, {\it The geometries of 3-manifolds},
 Bull. London Math. Soc. {\bf 15} (1983), 401--487.

 [So] T. Soma, Maps of non-zero degree to hyperboilc 3-manifolds. To appear
 J. Diff. Geom.

 [T] W.P. Thurston,  {\it The Geometry and Topology of
 3-manifolds}, Princeton University mimeo\-graphed notes (1979).

 [W] S.C. Wang, The existance of maps of non-zero degree between aspherical
 3-manifolds, Math. Zeit. {\bf 208} (1991), 147--160.

 [ZVC] H. Zieschang, E. Vogt and H. Coldeway, Surfaces and
 Planar Discontiuous Group, Lecture Note in Math. {\bf 835} Springer-Verlag
Berlin, 1980.

 \vskip 0.5 true cm

 \noindent Reid

 Department of Mathematics,

 University of Texas, Austin, TX 78712, USA

 \noindent Wang

 Department of Mathematics,

 Peking University, Beijing, 100871, P. R. China.

 \noindent Zhou

 Department of Mathematics,

 The East China Normal University,

 Shanghai, 200062, P.R. China

\bye